\documentclass{elsart}
\usepackage{stmaryrd}
\usepackage{amsfonts}
\usepackage{amsmath}
\usepackage{mathrsfs}
\usepackage{xcolor}
\usepackage[all]{xy}
\usepackage{graphicx}

\allowdisplaybreaks[4]
\newtheorem{tm}{Theorem}[section]
\newtheorem{pn}[tm]{Proposition}
\newtheorem{lm}[tm]{Lemma}

\theoremstyle{definition}
\newtheorem{dn}[tm]{Definition}
\newtheorem{rk}[tm]{Remark}
\newtheorem{ex}[tm]{Example}
\usepackage{enumerate}

\journal{Elsevier}
\usepackage{bm}

\begin{document}
	
\begin{frontmatter}
\title{Sober $L$-convex spaces and  $L$-join-semilattices}
\author{Guojun Wu, Wei Yao}
\address{School of Mathematics and Statistics, Nanjing University of Information Science and Technology, Nanjing, 210044, China\\
Applied Mathematics Center of Jiangsu Province, Nanjing University of Information Science and Technology, Nanjing, 210044, China\\
}

\address{}
\date{}

\begin{abstract}
With a complete residuated lattice $L$ as the truth value table, we  extend the definition of sobriety of classical convex spaces to the framework of $L$-convex spaces. We provide a specific construction for the sobrification of an $L$-convex space, demonstrating  that the full subcategory of sober $L$-convex spaces is  reflective in the category of $L$-convex spaces with  convexity-preserving mappings.  Additionally, we introduce the concept of Scott  $L$-convex structures on $L$-ordered sets. As an application of this type of sobriety, we obtain a characterization for the  $L$-join-semilattice completion of an $L$-ordered set:  an $L$-ordered set $Q$ is an $L$-join-semilattice completion of an $L$-ordered set $P$ if and only if the Scott  $L$-convex space  $(Q, \sigma^{\ast}(Q))$ is a sobrification of the Scott  $L$-convex space $(P, \sigma^{\ast}(P))$.
 \end{abstract}

\end{frontmatter}
\section{\bf Introduction }
A convex structure
 on a set  is a family closed under arbitrary intersections and directed unions, which  contains the empty
set as a member. Convex structure  can be
seen as the axiomatization of the usual convex sets in Euclidean
spaces.   Monograph \cite{Vanbook} provides an overview of the theory of  convex structures in detail.
Convex structure   exists in many branches of mathematics, including
lattices  \cite{Varlet,Van1984}, algebras \cite{Algebra,algebra4}, metric spaces \cite{Menger}, graphs \cite{graph,graph1,Soltan} and topological
spaces \cite{topology2,topology1}.

Analogous to topology, a convex structure is essentially  a family of subsets of a  background set. It is well-known that there are  close and rich connections between topological structures and ordered structures.  These connections are mainly reflected in domain theory \cite{Domain} and locale theory \cite{Stone}. In general, there are two primary  methods  to establish the connections  between orders and  topologies. { One approach  involves using intrinsic topologies on posets, such as the Scott topology, and the specialization orders induced by given topologies to link topological structures with ordered structures.  Famous results include the categorical isomorphism between injective $T_0$ spaces and continuous lattices \cite{Domain}.} Another one  takes  a pointless approach,  disregarding the underlying set. Specifically, we  equip the family of open sets   with set-inclusion order to  obtain a complete lattice; conversely, we can define an appropriate  spectral topology on a complete lattice. This method gives rise to category dualities between topological structures and ordered structures, such as the well-known duality   between sober topological spaces and spatial frames  \cite{Stone}.

Inspired by the above two approaches,  scholars have studied the intersection of convex structures and ordered structures. In \cite{Shen}, Shen et al. studied
the pointfree structure of convex spaces, demonstrating that  sober convex spaces and
algebraic lattices are categorically dual equivalent to each other. Recently, Xia \cite{Xia-further} extended these findings, discovering additional {categorical dualities} between convex structures and ordered structures in pointfree convex geometry.
The specialization order approach  has also been employed to study the ordered properties of convex structure.  In \cite{Jankowski}, Jankowski showed that  the category of injective $S_0$-convex spaces and that of frames are   isomorphic.   Yao and Zhou in \cite{Yao-Zhou}  established a categorical isomorphism
between sober convex spaces and join-semilattices by means with
 specialization order. These works   reveal  that there are close relationships between ordered structures and convex structures.

With the development of fuzzy mathematics,  Shen and  Shi in \cite{ShenShi} highlighted  that the fuzzy extensions of convex spaces are significant for both theoretical research and practical applications. In 1994, the notions of fuzzy convex spaces and hull operators were
first proposed by Rosa \cite{RosaDc}, which are called $I$-convex $(I = [0, 1])$ structures nowadays. Later, Maruyama \cite{Maruyama} extended Rosa's definition to
completely distributive lattice-valued setting and obtained the notion of $L$-convex
spaces. In recent years, Shi, Pang, and their team members have made substantial contributions to the study of $L$-convex structures. In \cite{PangShi}, Pang and Shi introduced several kinds of $L$-convex spaces
and discussed their categorical relationships.  In  \cite{ShenShi}, Shen and Shi  introduced the notions of $L$-convex systems and Scott-hull spaces,   establishing a categorical isomorphism between them. Similar to the bases and subbases of $L$-topology, Pang and Xiu in \cite{Pang3} studied the bases and subbases of $L$-convex structures. For further studies on fuzzy convex spaces, see  \cite{Pang2,Xiu2,WangKai}.

Just as the the intersection of orders and topologies can be extended to the fuzzy setting (see \cite{YaoTFS,PartII,Yao-Yue-FFS,Zhang-Shi-Li}), the  intersection of ordered structures and convex structures  can also be extended to the fuzzy setting. Yao and Zhou in \cite{Yao-Zhou-2021}  introduced sobriety of quantale-valued fuzzy convex spaces and
algebraicness of fuzzy complete lattices and then they established a  dually
equivalent  between the  category
of sober fuzzy convex spaces and that of algebraic fuzzy complete lattices, which is an extension of the work in  \cite{Shen}.
Motivated by {Yao's} work \cite{YaoTFS}, Xia in \cite{Xia-FFS} extended the Jankowski's work   \cite{Jankowski}  and established a categorical isomorphism   between the category of balanced $L$-$S_0$-convex spaces
 and that of fuzzy frames. Following  Xia' work, Liu, Yue and Wei in \cite{Liu-Meng-Yue} studied frame-valued Scott open set monad and proved that the related Kleisli monoids are precisely  the strong frame-valued convex spaces. Building on this, Pang in \cite{Pang-2023}   further developed the monadic approach in the theory of fuzzy convex structures.

{Recently, Liu and Yue  \cite{Liu-Yue-2024}  introduced algebraic irreducible convex sets (called compact  convex sets in this paper) with the
help of the fuzzy inclusion order between the convex sets in a $L$-convex space. This notion allowed them to  extended the theory of sober convex spaces to the fuzzy setting. In the classical case, the notions of  polytopes and compact convex sets are equivalent.
 It is thus natural to ask {\em whether one can postulate polytopes  in a $L$-convex spaces and then introduce another type of  sobriety of $L$-convex space with the help  of polytopes instead of   compact  convex sets}.
 Fortunately, this paper provides an affirmative answer. Moreover, unlike  Liu and Yue's sobriety, {\em our sobriety aims to foster deeper connections  between fuzzy ordered structures and fuzzy convex structures}, inspired by the role of topological sobriety  played in order theory \cite{XU-TopAPpp}.}

 In this paper,
we use  a complete residuated lattice $L$ as the truth value table.  This paper is organized as follows:  In Section 2, we recall basic concepts and results about lattices, $L$-orders and $L$-convex spaces.  In Section 3, we propose  a new type of sobriety of $L$-convex spaces and  provide a specific construction for the sobrification of any given $L$-convex space.  In Section 4,   we introduce the notion of $L$-join-semilattices and obtain     a construction for the $L$-join-semilattice completion  of any given $L$-ordered set via sobrification.

\section{\bf Preliminaries}
We refer to \cite{residuated,Hajek,Quantale} for contents on residuated lattices. We refer to \cite{MVTop} for notions of fuzzy sets, and to \cite{YaoTFS,PartI,Yao-Yue-FFS,FanZhang} for contents of fuzzy posets.

We say a subset $D$ of a poset $P$ is {\em directed} provided
it is non-empty and every two elements of $D$ has an upper bound in $D$. For $A\subseteq P$, write $\bigvee A$ for  the least upper
bound of $A$ and $\bigwedge A$ for the greatest lower bound of $A$.  In particular, we use the
convenient notation $x =\bigvee^{\uparrow}D$ to denote that the set $D$ is directed and $x$ is its
least upper bound.

Let $L$ be a complete lattice with  a bottom element $0$ and a top element 1 and let $\otimes$ be a binary operation
on $L$ such that $(L,\otimes, 1)$ is a commutative monoid.  The pair $(L,\otimes)$ is called a {\it complete residuated lattice}, {or a {\em commutative and integral quantale}},
if the operation $\otimes$ is distributive over joins; that is,
$$a\otimes(\bigvee S)=\bigvee_{s\in S} (a\otimes s).$$
For a complete residuated lattice $(L,\otimes)$, the operation
$\otimes$ gives rises to a right adjoint $\rightarrow:L\times L\longrightarrow L$
via the adjoint property  $$a\otimes b\leq c\ \Longleftrightarrow\ a\leq b\rightarrow c\ (\forall
a,b,c\in L).$$

%

\begin{lm}\label{lm-resi-lat}{\rm (\cite{residuated,Hajek})}Suppose that $(L,\otimes)$ is a complete residuated lattice. Then for all $a,b,c\in L,\ \{a_i|\ i\in I\},\ \{b_j|\ j\in J\}\subseteq L$,

{\rm (1)} $1= a\rightarrow b\Longleftrightarrow a\leq
b$;

{\rm (2)} $1\rightarrow a=a$;

{\rm (3)} $a\otimes(a\rightarrow b)\leq b$;
%
%

{\rm (4)} $a\rightarrow (b\rightarrow c) = (a\otimes b)\rightarrow c$

{\rm (5)} $(\bigvee_{i\in I}a_i)\rightarrow b=\bigwedge_{i\in I}(a_i\rightarrow
b)$;

{\rm (6)} $a\rightarrow(\bigwedge_{j\in J}b_j)=\bigwedge_{j\in J}(a\rightarrow
b_j)$.
%


\end{lm}

Every mapping $A:X\longrightarrow L$ is called an $L$-{\it subset} of $X$ and we  use $L^X$ to denote the collection of  $L$-subsets of $X$. Customarily, the crisp order $\leq$ on $L^X$ is  defined  pointwisely; that is $A\leq B\Leftrightarrow A(x)\leq B(x)\ (\forall x\in X)$. An $L$-subset $A$ is said to be {\it nonempty} if $\bigvee_{x\in X}A(x)=1$. Let $Y\subseteq X$ and $A\in L^X$, define $A|_{Y}\in L^Y$ by $A|_{Y}(y)=A(y)$ $(\forall y\in Y)$.
For an element $a\in L$, the notation $a_X$ denotes the constant $L$-subset of $X$ with  the value $a$, i.e., $a_X(x)=a\ (\forall x\in X)$.
For all $a\in L$ and $A\in  L^{X}$,  write $a\otimes A$ , for the $L$-subset given by $(a\otimes A)(x) = a\otimes A(x)$.

For each $a\in X$ and $Z\subseteq X\subseteq X^{\prime}$. Define characteristic functions  $1_a$ , $\chi_{Z}\in L^X$   by

\centerline{
$
\ \ 1_a(x)=\left\{\begin{array}{ll}1,& x=a;\\
0,&  x\neq a,
\end{array}\right.
$
$
\ \ \chi_{Z}(x)=\left\{\begin{array}{ll}1,& x\in Z;\\
0,&  x\notin Z.
\end{array}\right.
$
}
 It is worth noting that in this paper, we also use  {the symbols} $1_a$ , $\chi_{Z}\in L^{X^{\prime}}$ to denote the characteristic functions, { with the only difference being the domains of the mappings}. Readers should determine the domain of each characteristic function from the context to avoid any confusion.
\begin{dn}{(\cite{PartI,FanZhang})}
A mapping $e:P\times P\longrightarrow L$ is called
 an {\em $L$-order} if

{\rm (E1)} $\forall x\in    { P}$, $e(x,x)=1$;

{\rm (E2)} $\forall x,y,z\in { P}$, $e(x,y)\otimes e(y,z)\leq e(x,z)$;

{\rm (E3)} $\forall x,y\in { P}$, if $e(x,y)\wedge e(y,x)=1$,
then $x=y$.

The pair $(P,e)$ is called an $L$-{\em ordered set}. It is customary to write $P$ for the pair $(P, e)$.
\end{dn}
To avoid confusion, we sometimes use $e_p$ to denote the {$L$-order} on $P$.
A mapping  $f: P\longrightarrow Q$  between two $L$-ordered sets is said to be {\it $L$-order-preserving}  if for all $x, y\in P$, $e_{P}(x, y)\leq e_{Q}(f(x), f(y))$; $f$  is said to be  {\it $L$-order-isomorphic} if $f$  { is a bijection} and for all $x, y\in P$, $ e_{P}(x, y)= e_{Q}(f(x), f(y))$.

\begin{ex}\label{ex-L-ord}{(\cite{YaoTFS})}

{\rm(1)} Define $e_L:L\times L\longrightarrow L$ by $e_L(x, y)= x\to y$ $(\forall x, y\in L)$.   Then $e_L$ is an $L$-order on $L$.

{\rm(2)}  Define ${\rm sub}_X :L^X\times L^X\longrightarrow L$ by
 $${\rm sub}_X(A,B)=\bigwedge\limits_{x\in X}A(x)\rightarrow B(x)\ (\forall A, B\in L^X).$$ Then ${\rm sub}_X$ is an $L$-order on $L^X$, which is called the  {\em   inclusion $L$-order} on $L^X$.
If the background set is clear, then we always drop the subscript $X$ to be ${\rm sub}$.
\end{ex}

Let $f:X\longrightarrow Y$ be a mapping between two sets. The {\it Zadeh extensions} $f^\rightarrow:L^X\longrightarrow L^Y$ and $f^\leftarrow:L^Y\longrightarrow L^X$ are respectively given by
\vskip 4pt
\centerline{$f^\rightarrow(A)(y)=\bigvee\limits_{f(x)=y}A(x)\ (\forall A\in L^X)$,\   \ \qquad $f^\leftarrow(B)=B\circ f\ (\forall B\in L^Y).$}
\vskip 5pt

\begin{lm}{\em(\cite{DX-Zhang-2018})} For each mapping $f: X\longrightarrow Y$,

$(1)$ $f^{\rightarrow}: (L^{X}, {\rm sub}_X)\longrightarrow(L^{Y}, {\rm sub}_Y)$ is $L$-order-preserving;

$(2)$ $f^{\leftarrow}: (L^{Y}, {\rm sub}_Y)\longrightarrow(L^{X}, {\rm sub}_X)$ is $L$-order-preserving;

$(3)$ $f^{\rightarrow}$ is left adjoint to $f^{\leftarrow}$, written $f^{\rightarrow}\dashv f^{\leftarrow}$, that is
$$ {\rm sub}_Y(f^{\rightarrow}(A), B)={\rm sub}_X(A, f^{\leftarrow}(B))\ (\forall A\in L^X, B\in L^Y). $$
\end{lm}

Define ${\uparrow}x$ and ${\downarrow}x$ respectively by ${\uparrow}x(y)=e(x,y)$, ${\downarrow}x(y)=e(y,x)\ (\forall x,y\in P)$. An $L$-subset $S\in L^P$ is called a {\it lower set} (resp., an {\it upper set}) if $S(x)\otimes e(y,x)\leq S(y)$ (resp., $S(x){\otimes} e(x,y)\leq S(y)$) for all $x,y\in P$. Clearly,  ${\downarrow}x$ (resp.,  ${\uparrow}x$) is a lower (resp., an upper) set for every  $x\in P$.
\begin{dn}{(\cite{PartI,FanZhang})} Let $P$ be an $L$-ordered set. An element $x\in P$ is called a {\it supremum} of $A\in L^P$, denoted by $x=\sqcup A$, if
$$e(x,y)={\rm sub}(A,{\downarrow}y)\ (\forall y\in P).$$
Dually, an element $x$ is called an {\it infimum} of $A\in L^P$, denoted by $x=\sqcap A$, if
$$e(y,x)={\rm sub}(A,{\uparrow}y)\ (\forall y\in P).$$
\end{dn}

It is easy to check that if the supremum (resp., infimum) of an $L$-subset in an $L$-ordered set exists, then it must be unique.

For basic contents of $L$-convex spaces, we refer to \cite{PangShi,ShenShi}.
\begin{dn}\label{dn-con}{\rm (\cite{Liu-Yue-2024})} Let $X$ be a set and $\mathcal{C}\subseteq L^X$. The family $\mathcal{C}$ is called an \emph{ $L$-convex structure} on $X$ if it satisfies the following conditions:

{\rm (C1)} $0_X, 1_X\in \mathcal{C}$;

{\rm (C2)} $\bigvee_{i\in I}^{\uparrow}C_i\in \mathcal{C}$ for every directed subset  $\{ C_i\mid i\in I\}$  of $\mathcal{C}$;

{\rm (C3)}  $\bigwedge_{j\in J} C_j\in \mathcal{C}$ for every subset $\{ C_j\mid j\in J\}$ of $\mathcal{C}$;

{\rm (C4)}  $a\to C\in \mathcal{C}$ for all $a\in L$ and $C\in \mathcal{C}$.

The pair $(X, \mathcal{C})$ is called an {\em $L$-convex space}; every element of $\mathcal{C}$ is called a {\em convex set} of $(X, \mathcal{C})$.
\end{dn}
The standard name of  the $L$-convex space defined above is  {\em stratified
$L$-convex space} in the sense of   \cite{Liu-Yue-2024}. While in this paper, every $L$-convex space is always assumed
to be stratified,  so we omit the word ``stratified''.
When no confusion can arise, we often write
$X$ instead of $(X,\mathcal{C})$ for an $L$-convex space and write $\mathcal{C}(X)$ for the $L$-convex structure of $X$.

\begin{dn}{\rm (\cite{ShenShi})} Let $X$ be an $L$-convex space. Define a mapping $co_X: L^X\longrightarrow L^X$ by
$$ co_X(A)=\bigwedge\{B\in \mathcal{C}(X)\mid A\leq B\}\ (\forall A\in L^X),$$
called the {\em hull operator} of $(X,\mathcal{C}(X))$. For simplicity of notation, we always write $co$ instead of $co_X$ when no confusion can arise.
\end{dn}

%
%
%
%
%

 \begin{lm}{\rm (\cite[Proposition 2.10]{Liu-Yue-2024})}\label{pn-co} Let $X$ be an $L$-convex space. Then

 $(1)$ $a\otimes co(A)\leq co(a\otimes A)\ (\forall a\in L, A\in L^X)$;

 $(2)$ the hull operator $co:(L^X, {\rm sub})\longrightarrow (L^X, {\rm sub})$ is $L$-order-preserving;

 $(3)$ ${\rm sub}(A, B)={\rm sub}(co(A), B)\ (\forall A\in L^X, B\in \mathcal{C}(X))$.

\end{lm}
%
%
%

\begin{dn}{\rm (\cite{PangShi})}  Let $f: X\longrightarrow Y$ be a mapping between two {$L$-convex spaces}.
Then {$f$} is called

{\rm (1)} {\em convexity-preserving}  if for every $B \in \mathcal{C}(Y)$, $f^{ \leftarrow}(B)\in \mathcal{C}(X)$;

{\rm (2)} {\em convex-to-convex}  if for every $A\in \mathcal{C}(X)$, $f^{ \rightarrow}(A)\in \mathcal{C}(Y)$;

{\rm (3)} {\em convex-homeomorphic} if it is bijective, convexity-preserving and convex-to-convex.
\end{dn}

We say that  $X$ is {{\em convex-homeomorphic}}  to $Y$ if there exists {{\em a convex-homeomorphism}}
between $X$ and $Y$.

\begin{lm}{\rm (\cite{PangShi})}
Let $f:X\longrightarrow Y$ be a mapping between two $L$-convex spaces. Then $f$ is convexity-preserving if and only if
{$f^{\rightarrow}(co_X(A))\subseteq co_Y(f^{\rightarrow}(A))$ for every $A\in L^X$.}

\end{lm}

\section{Sober $L$-convex spaces}
Let $(X,\mathcal{C})$ be a classical convex space. A subset $E$ is called a {\it polytope} if there exists a nonempty finite subset $F\subseteq X$ such that $E =co(F)$. The space $X$  is called {\it sober} if, for every polytope $E$, there is a unique element $x$ such that $F=co(x)$. In order to extend the theory of sober convex spaces to the fuzzy setting, the first step is to postulate polytopes of $L$-convex spaces. Note  that a   subset $F$ of $X$ is finite if and only if, for every directed family $\{A_i\mid i\in I\}$, $F\subseteq\bigcup_{i\in I}^{\uparrow}A_i  $ implies that there exists $i\in I$ such that $F\subseteq A_i$.
Fortunately, making use of the fuzzy inclusion order between  $L$-subsets, we will naturally define polytopes in $L$-convex spaces and establish a theory of sober $L$-convex spaces.


\begin{dn}\label{dn-plo-sob} $(1)$ An $L$-subset  $F\in L^X$ is  said to be {\em finite} if
$${\rm sub}(F, \bigvee\nolimits^{\uparrow}_{i\in I} A_{i})=\bigvee\nolimits^{\uparrow}_{i\in I}{\rm sub}(F, A_{i})$$
for every directed family $\{A_i\mid i\in I\}\subseteq L^X$. A convex set $C$ of $X$ is called a {\em polytope}
if it is the hull of a nonempty finite $L$-subset.

$(2)$ An $L$-convex space is said to be   {\em sober} if for every nonempty finite set $F\in L^X$,  there exists a unique $x\in X$ such that $co(F)=co(1_x)$.
\end{dn}

\begin{rk}\label{rk-fin-plo}{ When $L$ is a frame,  we assert that  for each crisp finite subset $F$ of $X$,  the  characteristic function $\chi_F$  is a  finite $L$-subset of $X$. In fact, for every directed family $\{A_i\mid i\in I\}\subseteq L^X$, due to the finiteness of $F$, we have
\begin{align*}
{\rm sub}(\chi_F, \bigvee\nolimits^{\uparrow}_{i\in I} A_{i})&=\bigwedge_{x\in F}\bigvee\nolimits^{\uparrow}_{i\in I} A_{i}(x)\\
&=\bigvee\nolimits^{\uparrow}_{i\in I} \bigwedge_{x\in F}A_{i}(x)\\
&=\bigvee\nolimits^{\uparrow}_{i\in I}{\rm sub}(\chi_F,  A_{i}).
\end{align*}
But conversely, a finite $L$-subset may not necessarily be a characteristic function of a crisp finite subset. For example, when both the background set $X$ and the truth value table $L$ are crisp finite sets, every $L$-subset of  $X$ is a finite $L$-subset of $X$.}
\end{rk}


\begin{dn}(\cite{Liu-Yue-2024}) Let $X$ be an $L$-convex space. A convex set  $K$ is  said to be \emph{compact} (called {\em algebraic irreducible} in \cite{Liu-Yue-2024}) if $\bigvee_{x\in X}K(x)=1$ and
$${\rm sub}(K, \bigvee\nolimits^{\uparrow}_{i\in I}C_{i})=\bigvee\nolimits^{\uparrow}_{i\in I}{\rm sub}(K, C_{i})$$
for every directed family $\{C_i\mid i\in I\}\subseteq\mathcal{C}(X)$.

\end{dn}
   { In classical setting,  monograph \cite{Vanbook}  shows that compact convex sets are equivalent to polytopes. In the fuzzy setting, their relationship is much more  complex.  By Lemma \ref{pn-co}(3), polytopes are clearly compact convex sets. However, it remains unclear whether the reverse holds; that is, whether all compact convex sets are polytopes. We leave it as an open question.}



\begin{ex} Let $L=([0, 1], \otimes)$ be  a complete residuated lattice with $\otimes$ being $\wedge$. { Now, $L$ is a frame.}  Define a {stratified} $L$-convex structure
$$\mathcal{C}=\{a\wedge \phi\mid   a\in [0, 1], \phi: [0, 1]\longrightarrow [0, 1] \text{ is increasing}, \phi\geq id                 \}$$
 on $[0, 1]$.
{Specifically, a function $\mu: [0, 1]\to [0, 1]$ is a member of  $ \mathcal{C}$ if and only if $\mu$ is an increasing function; and there exists some $a\in [0,1]$ such that
$\mu(x)\geq x$ if $x\in [0, a)$ and $\mu(x)=a$ if $x\in [a, 1]$.}

We next show that $A$ is nonempty finite $L$-subset of $[0, 1]$ if and only if   there exists a nonempty finite subset $F\subseteq_{fin} [0,1]$ such that $A =\chi_{F}$. By Remark \ref{rk-fin-plo}, it remains to prove the ``only if" part. We  divided this proof into two steps.

{\bf Step 1.} We prove that $\{x\in [0,1]\mid A(x)\neq 0\}$ is a nonempty finite set.
Since $\bigvee_{x\in X}A(x)=1$, there exists $x\in [0,1]$ such that $A(x)\neq 0$.  Assume that $\{x\in [0,1]\mid A(x)\neq 0\}$ is an infinite set. Then for every $F\subseteq_{fin} X$, there exists $x$ such that $A(x)\neq 0$, but $x\notin F$.  Noticing that $\otimes$ is $\wedge$, we have $\bigvee_{F\subseteq_{fin} [0,1]}{\rm sub}(A, A\wedge \chi_{F})=0$. Since $A$ is finite, it holds that
$$1={\rm sub}(A, \bigvee_{F\subseteq_{fin}[0,1]}A\wedge \chi_{F})= \bigvee_{F\subseteq_{fin} [0,1]}{\rm sub}(A, A\wedge \chi_{F}), $$ a contradiction. Thus, $\{x\in [0,1]\mid A(x)\neq 0\}$ is a nonempty finite set.

{\bf Step 2.} We prove that there is no $x_0$ such that $A(x_0)\in (0,1)$.

Assume that there exists $x_0$ such that $A(x_0)\in (0,1)$. Write $A(x_0)=y_0$. It is clear that there exists $t_0> 1$, such that $y_0-\frac{1}{t_0}> 0$. Define $A_n:[0, 1]\longrightarrow [0, 1]$ by  $A_n(x_0)=y_0-\frac{1}{t_0+n}$; when $x\neq x_0$, $A_n(x)=A(x)$. We obtain a directed family $ \{A_n\mid n\in \mathbb{N}  \}   $. Obviously, $\bigvee_{n\in \mathbb{N}}A_n=A$; that is to say ${\rm sub}(A, \bigvee_{n\in \mathbb{N}}A_n)=1$. For every $n$, it holds that
$${\rm sub}(A, A_n)=A(x_0)\rightarrow A_n(x_0)=y_0-\frac{1}{t_0+n}.$$
 Therefore, $\bigvee_{n\in \mathbb{N}}{\rm sub}(A, A_n)=y_0< 1$, a contradiction.

Thus, there exists a nonempty finite subset $F\subseteq_{fin} [0,1]$ such that $A =\chi_{F}$.

Write $min (F)=b_0$. It is routine to check that the smallest convex set containing $A$ is
$$
\ \ \phi(x)=\left\{\begin{array}{ll}1,&x\geq b_0;\\
x,& x< b_0.
\end{array}\right.
$$
It is clear that  $co(A)=co(1_{b_0})$. Thus $([0,1],\mathcal{ C})$ is a sober $L$-convex space. {By the way, it is mechanical to check that in $([0,1],\mathcal{ C})$, every compact convex  set is exactly a polytope; we leave this verification to the reader.}\hfill$\Box$\\


\end{ex}

In this paper, we use ${\rm cp}(\mathcal{C}(X))$  denote the set of all compact convex sets of $X$. For each $A\in \mathcal{C}(X)$, define
$$ \phi(A):{\rm cp}(\mathcal{C}(X))\longrightarrow L  $$
by
$\phi(A)(K)={\rm sub}(K, A)$.

\begin{lm}\label{lm-spec-conv} Let $X$ be an $L$-convex space.

$(1)$ $\phi(a_X)(K)=a$ for $a=0, 1$ and for every $K\in {\rm cp}(\mathcal{C}(X))$;

$(2)$ $\phi(\bigvee_{i\in I}^{\uparrow}C_i)=\bigvee_{i\in I}^{\uparrow}\phi(C_i)$ for every directed family $\{C_i\mid i\in I\}\subseteq  \mathcal{C}(X)$ ;

$(3)$ $\phi(\bigwedge_{j\in J}C_j)=\bigwedge_{j\in J}\phi(C_j)$ for every  family $\{C_j\mid j\in J\}\subseteq  \mathcal{C}(X)$ ;

$(4)$ $\phi(a\to C)=a\to \phi(C)$ for $a\in L$ and $C\in \mathcal{C}(X)$;

$(5)$ ${\rm sub}_X(A, B)={\rm sub}_{{\rm cp}(\mathcal{C}(X))}(\phi(A), \phi(B))$ for all $A, B\in \mathcal{C}(X)$.
\end{lm}
\noindent {\bf Proof.}
%
%
%
The verification is straightforward  by Lemma \ref{lm-resi-lat}.\hfill$\Box$\\

Given an $L$-convex space $X$, by Lemma \ref{lm-spec-conv}, it is straightforward to check that
$\{\phi(A)\mid A\in \mathcal{C}(X)\}$
 is an $L$-convex space on ${\rm cp}(\mathcal{C}(X))$ and denoted by $\mathcal{C}({\rm cp}(\mathcal{C}(X)))$. We write ${\rm Cp}(\mathcal{C}(X))$, rather than $({\rm cp}(\mathcal{C}(X)), \mathcal{C}({\rm cp}(\mathcal{C}(X))))$, for the resulting $L$-convex space. By Lemma \ref{lm-spec-conv}(2)(5), it is easy to check that $A$ is a compact convex set of $X$ if and only if $\phi(A)$ is a compact convex set of {${\rm Cp}(\mathcal{C}(X))$}.

\begin{pn}\label{pn-cp-sob} Let $X$ be an $L$-convex space. Then ${\rm Cp}(\mathcal{C}(X))$ is a sober $L$-convex space.
\end{pn}
\noindent {\bf Proof.} We first show that $\phi(K)=co(1_{K})$ for every $K\in {\rm cp}(\mathcal{C}(X))$. Since $\phi(K)(K)=1$, we have $1_K\leq\phi(K)$. Thus $co(1_K)\leq\phi(K)$. For each convex set $A$, if $1_K\leq \phi(A)$, then $\phi(A)(K)={\rm sub}(K,A)=1$. For every $G\in {\rm cp}(\mathcal{C}(X))$,
$$\phi(K)(G)={\rm sub}(G,K)={\rm sub}(G, K)\wedge {\rm sub}(K, A)\leq {\rm sub}(G,A)=\phi(A)(G).  $$
Thus $\phi(K)\leq\phi(A)$. This shows that $\phi(K)=co(1_{K})$.

Let $\phi(A)\in \mathcal{C}({\rm cp}(\mathcal{C}(X)))$ be a polytope. Then $\phi(A)$ is a compact convex set of  ${\rm Cp}(\mathcal{C}(X))$ and $A$ is a compact convex set of  $X$.
 Therefore, $\phi(A)=co(1_{A})$. The  uniqueness of $A$ can be obtained directly from Lemma \ref{lm-spec-conv}(5).  Thus ${\rm Cp}(\mathcal{C}(X))$ is a sober $L$-convex space.\hfill$\Box$\\
In the following second part of this section, we will provide a construction of  the sobrification of an $L$-convex space. We first need a definition as follows.

\begin{dn}\label{dn-F-close} Let $X$ be an $L$-convex space and $A\subseteq X$. Then $A$ is called an {{\em F-closed set}} if for every nonempty finite $L$-subset with $F\leq \chi_A$ and $co_X(F)=co_X(1_x)$ implies $x\in A$.
\end{dn}

In this paper, let $\mathfrak{F}(X)$ denote the family of F-closed sets of $X$.  The collection $\mathfrak{F}(X)$ is a closure system on $X$, that is to say $\mathfrak{F}(X)$ is closed under arbitrary intersections (including empty intersection). The collection $\mathfrak{F}(X)$ gives rise to a closure operator ${\rm cl}_{\mathfrak{F}}:\mathcal{P}(X)\longrightarrow \mathcal{P}(X)$ defined by
$${\rm cl}_{\mathfrak{F}}(B)=\bigcap\{A\in \mathfrak{F}(X)\mid B\subseteq A\}.$$
We call ${\rm cl}_{\mathfrak{F}}(B)$ the {\it F-closure} of $B$.

Let $f:X\longrightarrow Y$ be a mapping between two  $L$-convex spaces. Then $f$ is said to be {\it F-continuous} if for every $A\in \mathfrak{F}(Y)$, one has $f^{-1}(A)\in \mathfrak{F}(X)$.

Recall that an $L$-convex space $X$ is said to be $S_0$ if for all $x, y\in X$, $co(1_x)=co(1_y)$ implies  $x=y$.
\begin{pn}\label{pn-fcon} Let $X$ and $Y$ be two   $L$-convex spaces and let $f:X\longrightarrow Y$ be convexity-preserving. The following statements hold:

$(1)$ for every nonempty $L$-subset $A\in L^X$, $co_Y(f^{\rightarrow}(A))=co_Y(1_{f(x)})$, where $co_X(A)=co_X(1_x)$;

$(2)$ $f$ is F-continuous;

$(3)$ if $Y$ is $S_0$, $g: X\longrightarrow Y$ is  convexity-preserving and $Z\subseteq X$ with $g|_Z=f|_Z$, then $g|_{{\rm cl}_{\mathfrak{F}}(Z)}=f|_{{\rm cl}_{\mathfrak{F}}(Z)}$.
\end{pn}
\noindent {\bf Proof.}
(1) Since $f:X\longrightarrow Y$ be convexity-preserving, we have
\begin{align*}
1=&{\rm sub}(f^{\rightarrow}(A), f^{\rightarrow}(co_X(A)))\\
=& {\rm sub}(f^{\rightarrow}(A), f^{\rightarrow}(co_X(1_x)))\\
\leq& {\rm sub}(f^{\rightarrow}(A), co_Y(f^{\rightarrow}(1_x)))\\
=&{\rm sub}(co_Y(f^{\rightarrow}(A)), co_Y(f^{\rightarrow}(1_x)))\\
=&{\rm sub}(co_Y(f^{\rightarrow}(A)), co_Y(1_{f(x)})
\end{align*}
Thus $co_Y(f^{\rightarrow}(A))\leq co_Y(1_{f(x)})$.

On the other hand, since
$$1=f^{\rightarrow}(co_X(1_x))(f(x))=f^{\rightarrow}(co_X(A))(f(x))\leq  co_Y(f^{\rightarrow}(A))(f(x)), $$
we have $co_Y(1_{f(x)})\leq co_Y(f^{\rightarrow}(A))$. Thus $co_Y(f^{\rightarrow}(A))=co_Y(1_{f(x)})$.

(2) Let $B$ be an F-closed set of $Y$. We will show that $f^{-1}(B)$ is an F-closed set of $X$. For nonempty set $F$ with $F\leq \chi_{f^{-1}(B)}$ and $co_X(F)=co_X(1_x)$, we have
$$F\leq \chi_{f^{-1}(B)}=\chi_B\circ f=f^{\leftarrow}(\chi_B).$$
By $f^{\rightarrow}\dashv  f^{\leftarrow}$, we have $f^{\rightarrow}(F)\leq\chi_B$ and $f^{\rightarrow}(F)$ is a nonempty finite $L$-subset of $Y$.
By  Part (1), we have $co_Y(f^{\rightarrow}(F))= co_Y(1_{f(x)})$. Since $B$ is an F-closed set, we have $f(x)\in B$. Therefore, $x\in f^{-1}(B)$. This shows that $f^{-1}(B)$ is an F-closed set of $X$. Thus $f$ is F-continuous.

(3) Write $M=\{x\in X\mid g(x)=f(x)\}$. Clearly, $Z\subseteq M$. We will show that $M$ is an F-closed set of $X$. Let $F$ be a nonempty $L$-subset with $F\leq\chi_M$ and $co_X(F)=co_X(1_{x_0})$. By Part (1), we have
$$  co_Y(1_{f(x_{0})})= co_Y(f^{\rightarrow}(F))=co_Y(g^{\rightarrow}(F))=co_Y(1_{g(x_{0})}).$$
Since $Y$ is $S_0$, it follows that $f(x_{0})=g(x_{0})$, i.e., $x_0\in M$. This shows that $M$ is an F-closed set and ${\rm cl}_{\mathfrak{F}}(Z)\subseteq M $. Thus $g|_{{\rm cl}_{\mathfrak{F}}(Z)}=f|_{{\rm cl}_{\mathfrak{F}}(Z)}$.
 \hfill$\Box$
\vskip 3pt

Given an $L$-convex space $X$, define
$$\Theta(X)=\{co(1_x)\mid x\in X\}.$$ We use the symbol $X^F$ to denote the F-closure of $\Theta(X)$ in ${\rm Cp}(\mathcal{C}(X))$, i.e., $X^F:={\rm cl}_{\mathfrak{F}}(\Theta(X)) $. For $A\in \mathcal{C}(X)$, define $\varphi(A): X^F\longrightarrow L$ by $\varphi(A)(C)={\rm sub}(C,A)$, i.e., {$\varphi(A)=\phi(A)|_{X^F}$}.
It is easy to see  that $\{\varphi(A)\mid A\in \mathcal{C}(X)\}$ is an $L$-convex structure on $X^F$,  {denoted by $\mathcal{C}(X^F)$}. When we see $X^F$ as an $L$-convex space, the related $L$-convex structure is always assumed  to be the resulting $L$-convex structure. It is easy to see that $X^F$ is a subspace of ${\rm Cp}(\mathcal{C}(X))$.
Similar to Lemma \ref{lm-spec-conv}(5),  we also have  ${\rm sub}_X(A, B)={\rm sub}_{X^F}(\varphi(A), \varphi(B))$ for all $A, B\in \mathcal{C}(X)$, which will be useful.
%
%
%
%
%
%
%

\begin{pn} $X^F$ is a  sober $L$-convex space.
\end{pn}
\noindent {\bf Proof.}
Let $\mathcal{K}$ be a nonempty finite $L$-subset of $X^F$. Define $\mathcal{K}^{\prime}:{\rm cp}(\mathcal{C}(X))\longrightarrow L$ by
$$
\ \ \mathcal{K}^{\prime}(A)=\left\{\begin{array}{ll}\mathcal{K}(A),& A\in X^F;\\
\  \ 0,& A\in {\rm cp}(\mathcal{C}(X))-X^F.
\end{array}\right.
$$
We will show that $\mathcal{K}^{\prime}$ is a nonempty {finite} $L$-subset of ${\rm cp}(\mathcal{C}(X))$. It is obvious that  $\mathcal{K}^{\prime}$ is nonempty. For every directed family $\{\mathcal{A}_i\mid i\in I\}\subseteq L^{{\rm cp}(\mathcal{C}(X))}$, write $\mathcal{A}_i|_{X^F}=\mathcal{A}_i^{\ast}$. Then  $\{\mathcal{A}_i^{\ast}\mid i\in I\}$ is a directed  family of $L^{X^F}$. Then we have
\begin{align*}
{\rm sub}_{ {\rm cp}(\mathcal{C}(X))}(\mathcal{K}^{\prime}, \bigvee_{i\in I}\mathcal{A}_i)&= \bigwedge_{A\in {\rm cp}(\mathcal{C}(X))}\mathcal{K}^{\prime}(A)\to \bigvee_{i\in I}\mathcal{A}_i(A)\\
&= \bigwedge_{A\in X^F}\mathcal{K}(A)\to \bigvee_{i\in I}\mathcal{A}_i^{\ast}(A)\\
&={\rm sub}_{X^F}(\mathcal{K}, \bigvee_{i\in I}\mathcal{A}_i^{\ast})\\
&=\bigvee_{i\in I}{\rm sub}_{X^F}(\mathcal{K},\mathcal{A}_i^{\ast})\\
&=\bigvee_{i\in I}\bigwedge_{A\in X^F}\mathcal{K}(A)\to\mathcal{A}_i^{\ast}(A)\\
&= \bigvee_{i\in I}\bigwedge_{A\in {\rm cp}(\mathcal{C}(X))}\mathcal{K}^{\prime}(A)\to\mathcal{A}_i(A)\\
&=\bigvee_{i\in I}{\rm sub}_{{\rm cp}(\mathcal{C}(X))}(\mathcal{K}^{\prime}, \mathcal{A}_i).
\end{align*}
This shows that $\mathcal{K}^{\prime}$ is a nonempty {finite} $L$-subset of ${\rm cp}(\mathcal{C}(X))$. Therefore, $co_{{\rm cp}(\mathcal{C}(X))}(\mathcal{K}^{\prime})$ is a polytope in ${\rm Cp}(\mathcal{C}(X))$. It follows from the sobriety of ${\rm Cp}(\mathcal{C}(X))$  that  there exists a unique $A_0\in {\rm cp}(\mathcal{C}(X))$ such that
$$co_{{\rm cp}(\mathcal{C}(X))}(\mathcal{K}^{\prime})=co_{{\rm cp}(\mathcal{C}(X))}(1_{A_{0}})(=\phi(A_0)).$$
Since $X^F$ is an F-closed set in ${\rm Cp}(\mathcal{C}(X))$ and $\mathcal{K}^{\prime}\leq \chi_{X^F}\in L^{{\rm cp}(\mathcal{C}(X))}$, we have $A_{0}\in X^F$.
We claim that $co_{X^F}(\mathcal{K})=co_{{\rm cp}(\mathcal{C}(X))}(\mathcal{K}^{\prime})|_{X^F}$. In fact,
\begin{align*}
co_{X^F}(\mathcal{K})&= \bigwedge\{\varphi(C)\mid C\in \mathcal{C}, \mathcal{K}\leq \varphi(C)\}\\
&= \bigwedge\{\phi(C)|_{X^F}\mid C\in \mathcal{C}, \mathcal{K}^{\prime}\leq \phi(C)\}\\
&=(\bigwedge\{\phi(C)\mid C\in \mathcal{C}, \mathcal{K}^{\prime}\leq \phi(C)\})|_{X^F}\\
&=co_{{\rm cp}(\mathcal{C}(X))}(\mathcal{K}^{\prime})|_{X^F}.
\end{align*}
By this fact and $A_{0}\in X^F$,  it follows that
$$co_{X^F}(\mathcal{K})=co_{{\rm cp}(\mathcal{C}(X))}(\mathcal{K}^{\prime})|_{X^F}=\phi(A_0)|_{X^F}=\varphi(A_0). $$
Hence $ co_{X^F}(\mathcal{K})=co_{X^F} (1_{A_0})$, where the uniqueness of $A_0$ can be
derived from Lemma \ref{lm-spec-conv}(5). Thus $X^F$ is a  sober $L$-convex space.
\hfill$\Box$\\

\begin{pn}\label{pn-sober-iso}Let $X$ be an $L$-convex space. Define $\xi_{X}:X\longrightarrow X^F$ by $\xi_{X}(x)=co(1_x)$. Then

$(1)$ $\xi_{X}$ is a convexity-preserving mapping;

$(2)$ $X$ is sober if and only if $\xi_{X}$  is {a convex-homeomorphism}.
\end{pn}
\noindent {\bf Proof.}
(1)  Let $A\in \mathcal{C}(X)$ and $x\in X$. Then
$$ \xi_{X}^{\leftarrow}(\varphi(A))(x)=\varphi(A)(co(1_x))={\rm sub}(co(1_x), A)={\rm sub}(1_x, A)=A(x).$$
Therefore, $ \xi_{X}^{\leftarrow}(\varphi(A))=A$. Thus $\xi_{X}$ is a convexity-preserving mapping.

(2)
 Suppose that $\xi_{X}$  is {a convex-homeomorphism}. Since $X^F$ is sober, it is clear that $X$ is sober.
 Conversely,
 let $X$ be a sober space. Then  $\xi_{X}$ is an injection. We claim that $\Theta(X)$ is F-closed in ${{\rm Cp}(\mathcal{C}(X))}$, i.e., ${\rm cl}_{\mathfrak{F}}(\Theta(X))=\Theta(X)(=X^F)$. Let $\mathcal{K}$ be a nonempty finite $L$-subset of ${{\rm cp}(\mathcal{C}(X))}$ and $\mathcal{K}\leq \chi_{ \Theta(X)}$. Define $K\in L^X$ by $K(x)=\mathcal{K}(co_X(1_x))$.  It is routine to check that $K$ is a nonempty finite $L$-subset of $X$ and $\xi_{X}^{\rightarrow}(K)=\mathcal{K}$. Since $X$ is sober, it follows that there exists a unique $x_0\in X$ such that $co_X(K)=co_X(1_{x_0})$.
It follows from Part (1) and  Proposition \ref{pn-fcon}(1) that
$$co_{{{\rm cp}(\mathcal{C}(X))}}(\mathcal{K})=co_{{{\rm cp}(\mathcal{C}(X))}}(1_{co_X(1_{x_0})})(=\phi(co_X(1_{x_0}))).$$
Notice that $co_{X}(1_{x_0})\in \Theta(X)$. Thus $\Theta(X)$ is F-closed in ${{\rm Cp}(\mathcal{C}(X))}$. Hence ${\rm cl}_{\mathfrak{F}}(\Theta(X))=\Theta(X)(=X^F)$.
It follows that  $\xi_{X}$ is a bijection. It is routine to check that $\xi_{X}^{\rightarrow}(A)(co_X(1_x))=A(x)=\varphi (A)(co_X(1_x))$ for all $A\in \mathcal{C}(X)$.
Hence $\xi_{X}^{\rightarrow}(A)=\varphi (A)$. Thus $\xi_{X}$ is a {convex-homeomorphism}.
\hfill$\Box$\\

\begin{lm}\label{lm-Fclo-two} If $Z\subseteq X^F$ is an F-closed set of $X^F$, then $Z$ is an F-closed set of ${{\rm Cp}(\mathcal{C}(X))}$.
\end{lm}
\noindent {\bf Proof.} Let $K$ be a nonempty finite  $L$-subset of ${{\rm cp}(\mathcal{C}(X))}$ and $K\leq\chi_{Z}\in L^{{\rm cp}(\mathcal{C}(X))}$. Since ${{\rm Cp}(\mathcal{C}(X))}$ is  sober, there exists a unique $A_0\in {{\rm cp}(\mathcal{C}(X))}$ such that
$$co_{{{\rm cp}(\mathcal{C}(X))}}(K)=\phi(A_0)(=co_{{{\rm cp}(\mathcal{C}(X))}}(1_{A_0})).$$
 Since $X^F$ is an F-closed set of ${{\rm Cp}(\mathcal{C}(X))}$ and $K\leq \chi_{X^F}\in L^{{\rm cp}(\mathcal{C}(X))}$, we have $A_0\in X^F$.
Define $K^{\ast}:X^F\longrightarrow L$ by
$K^{\ast}(A)=K(A)$, i.e., $K^{\ast}=K|_{X^F}$. It is routine to check that  $K^{\ast}$ is a nonempty finite $L$-subset of $X^F$. Since $X^F$ is a subspace of ${{\rm Cp}(\mathcal{C}(X))}$, it is routine to obtain that
$$co_{X^F}(K^{\ast})=\varphi(A_0)(=co_{X^F}(1_{A_0})).$$ Since $Z$ is an F-closed set of $X^F$  and $K^{\ast}\leq\chi_{Z}\in L^{X^F}$, it follows that $A_0\in Z$.
 This shows that $Z$ is an F-closed set of ${{\rm Cp}(\mathcal{C}(X))}$.
\hfill$\Box$\\



We give a standard definition of sobrification as follows.
\begin{dn}
Let $X$ be an  $L$-convex space, let $Y$  be a sober  $L$-convex space and  let $j:X\longrightarrow Y$ be a  convexity-preserving mapping. Then  $(Y, j)$, or  $Y$ is called a {{\em sobrification}}  of $X$ if for every sober $L$-convex space $Z$ and every  convexity-preserving mapping $f:X\longrightarrow Z$, there exists a unique    convexity-preserving mapping $\overline{f}: Y\longrightarrow Z$ such that $f=\overline{f}\circ j$.

\end{dn}
By the  universal property of sobrifications, it is easy to see that the sobrification of an  $L$-convex space  is unique up to { convex-homeomorphism}.
 Next, we present the main result of this section as follows.
\begin{tm}\label{tm-sob-ri} $X^F$ with mapping $\xi_X$ is a sobrification of $X$.
\end{tm}
\noindent {\bf Proof.} Let $Y$ be a sober $L$-convex space and let $f:X\longrightarrow Y$  be a convexity-preserving mapping. Define $g: {\rm cp}(\mathcal{C}(X))\longrightarrow {\rm cp}(\mathcal{C}(Y))$ by $g(K)=co_{Y}(f^{\rightarrow}(K))$. For every $x\in X$,
$$g( \xi_X(x))=g(co_X(1_x))=co_Y(f^{\rightarrow}(co_X(1_x)))=co_Y(1_{f(x)}), $$
 Thus $g(co_X(1_x))=co_Y(1_{f(x)})$, which shows that $g(\Theta(X))\subseteq \Theta(Y)$.
For every $A\in {\rm cp}(\mathcal{C}(X))$ and $C\in \mathcal{C}(Y)$,
\begin{align*}
g^{\leftarrow}(\phi(C))(A)&=\phi(C)( co_{Y}(f^{\rightarrow}(A))\\
&={\rm sub}(co_{Y}(f^{\rightarrow}(A)), C)\\
&={\rm sub}(f^{\rightarrow}(A), C)\\
&={\rm sub}(A, f^{\leftarrow}(C))\\
&=\phi(f^{\leftarrow}(C))(A).
\end{align*}
Thus $g^{\leftarrow}(\phi(C))=\phi(f^{\leftarrow}(C))$. This shows that $g$ is convexity-preserving.
It follows from Proposition \ref{pn-sober-iso}(2) that $Y^F=\Theta(Y)$ which is an F-closed set of ${\rm cp}(\mathcal{C}(Y))$. By Proposition \ref{pn-fcon}(2),   $g^{-1}( \Theta(Y))$ is an F-closed set of ${\rm cp}(\mathcal{C}(X))$.  Thus $$X^F={\rm cl}_{\mathfrak{F}}(\Theta(X))\subseteq g^{-1}( \Theta(Y)).$$
Therefore $g(X^F)\subseteq  \Theta(Y)$. Thus for every $K\in X^F$, there exists a unique $y\in Y$ such that $g(K)=co_Y(1_y)$.  Define
$\overline{f}(K)$ to be this $y$. For every convex set $B$ of $Y$ and every $K\in X^{F}$,
we have
\begin{align*}
\overline{f}^{\leftarrow}(B)(K)&=B(\overline{f}(K))
={\rm sub}(co_Y(1_{\overline{f}(K)}), B)\\
&={\rm sub}({g(K)}, B)\\
&={\rm sub}(co_Y(f^{\rightarrow}(K)), B)\\
&={\rm sub}(f^{\rightarrow}(K),B)\\
&={\rm sub}(K,f^{\leftarrow}(B))\\
&=\varphi(f^{\leftarrow}(B))(K).
\end{align*}
{Thus, $\overline{f}^{\leftarrow}(B)=\varphi(f^{\leftarrow}(B))\in \mathcal{C}(X^F)$.} This shows that $\overline{f}: X^F\longrightarrow Y$ is convexity-preserving. Notice that for every $x\in X$, $g(co_X(1_x))=co_Y(1_{f(x)})$.
Thus $\overline{f}\circ \xi_X=f$.

In order to  show the uniqueness of $\overline{f}$, let $h:X^F\longrightarrow Y$ such that $f=h\circ\xi_X$. Then $\overline{f}|_{\Theta(X)}=h|_{\Theta(X)}$.
By Lemma \ref{lm-Fclo-two}, we know that the F-closure of $\Theta(X)$ in $X^F$  coincides with  that in ${\rm cp}(\mathcal{C}(X))$; that is $X^F$. Then by Proposition \ref{pn-fcon}(3), we know that $\overline{f}=h$.
\hfill$\Box$\\

We denote the category of $L$-convex spaces with convexity-preserving mappings as morphisms by $L$-{\bf CS}.
Theorem \ref{tm-sob-ri} shows that the full subcategory of sober $L$-convex spaces is reflective in $L$-{\bf CS}. For a detailed discussion on category theory, we refer the reader to \cite{Category}.

\section{$L$-join-semilattice completion via sobrification}
In this section,   we  introduce the notion of $L$-join-semilattices and   Scott $L$-convex structures on  $L$-ordered sets. By means with specialization $L$-order and Scott $L$-convex structure, we will discuss the relationship  between sober $L$-convex spaces and $L$-join-semilattices. Finally, we  also present     a construction for the $L$-join-semilattice completion  of any given $L$-ordered set via sobrifications. {These results were not reflected in the framework of Liu and Yue's sobreity \cite{Liu-Yue-2024}. }\\

Given an $S_0$ $L$-convex space $X$, define $e_{\mathcal{C}(X)}: X\times X\longrightarrow L$ by
$$e_{\mathcal{C}(X)}(x, y)=\bigwedge_{A\in \mathcal{C}(X)}A(y)\to A(x).$$
It is easy to see that $e_{\mathcal{C}(X)}$ is an $L$-order on $X$, called the {\it specialization $L$-order} of the space $X$. Accordingly, we will write $\Omega(X)$
for the $L$-ordered set obtained by equipping $X$ with its specialization $L$-order. In the following, unless otherwise specified, the $L$-order of a given $L$-convex space refers to its specialization $L$-order.

\begin{pn}\label{pn-spec-co} Let $X$ be an $L$-convex space. Then $e_{\mathcal{C}(X)}(x, y)=co(1_y)(x)$ for all $x, y\in X$.
\end{pn}
\noindent {\bf Proof.}
 On one hand, $$e_{\mathcal{C}(X)}(x, y)\leq co(1_y)(y)\to co(1_y)(x)=1\to  co(1_y)(x)=co(1_y)(x).$$
On the other hand, we have
\begin{align*}
co(1_y)(x)\otimes A(y)&= co(1_y)(x)\otimes {\rm sub}(1_y, A)\\
&=co(1_y)(x)\otimes {\rm sub}(co(1_y), A)
\leq A(x).
\end{align*}
It follows that
$$co(1_y)(x)\leq \bigwedge_{A\in \mathcal{C}(X)} A(y)\to A(x)=e_{\mathcal{C}(X)}(x, y).$$
The proof is finished.
\hfill$\Box$\\
\begin{pn}\label{pn-xf-spe} Let $X$ be an $L$-convex space. Then for every $A, B\in {\rm cp}(\mathcal{C}(X))$, $e_{{\rm cp}(\mathcal{C}(X))}(A, B)={\rm sub}(A, B)$; for every $A, B\in X^F$, $e_{X^F}(A, B)={\rm sub}(A, B)$.
\end{pn}
\noindent {\bf Proof.} By the proof of Proposition \ref{pn-cp-sob} and Proposition \ref{pn-spec-co}.
$$e_{{\rm cp}(\mathcal{C}(X))}(A, B)=co_{{\rm cp}(\mathcal{C}(X))}(1_B)(A)=\phi(B)(A)={\rm sub}(A, B).$$
Similarly, for every $A, B\in X^F$, $e_{X^F}(A, B)={\rm sub}(A, B)$.
\hfill$\Box$\\

Yue, Yao, and Ho  introduced the notion of  Scott convex structures on a join-semilattice (see \cite[Example 5.1]{Yue-Yao-Ho}), which precisely corresponds to the collection of all ideals of the join-semilattice. This notion of Scott convex structures on a join-semilattice can be straightforwardly extended to a more general poset. We now define  the Scott $L$-convex structure on an $L$-ordered set.

\begin{dn} Let $P$ be an $L$-ordered set. An $L$-subset $A\in L^P$ is called a {\em Scott  convex set} if $A$ is a lower set and for every  nonempty finite $L$-subset $F$ with supremum exists,  it holds that
$${\rm sub}(F, A)\leq A(\sqcup F).$$
\end{dn}

{It is to observe that if we replace the nonempty finite $L$-subsets with  directed $L$-subsets in the above definition,  we precisely obtain the definition of  Scott closed sets (see \cite[Definition 5.1]{Zhang-Wang}) in
 $L$-cotopological spaces.}
Let $\sigma^{\ast}(P)$ denote the set of all Scott $L$-convex sets of $P$. It is straightforward to check that   $\sigma^{\ast}(P)$ is an actual $S_0$ $L$-convex structure and we call  $\sigma^{\ast}(P)$ the {\em Scott $L$-convex structure} of $P$. {It is easy to check that for an $L$-ordered set $(P, e)$,  $\Omega(P, \sigma^{\ast}(P))=(P, e)$.}

As a counterpart of  Scott continuous mappings between two $L$-ordered sets (see \cite[Definition 2.12]{PartII}), we give the following definition.
\begin{dn} An $L$-order-preserving mapping $f:P\longrightarrow Q$ is said to be {\em Scott convexity-preserving} if for every nonempty finite $L$-subset $F$ with a supremum, the supremum of $f^\rightarrow(F)$  exists and $f(\sqcup F)=\sqcup f^\rightarrow(F)$.
\end{dn}
The following shows that Scott convexity-preserving mappings consistent   with the convexity-preserving mappings between the related Scott $L$-convex spaces.
\begin{pn}\label{pn-sco-dir} $f:P\longrightarrow Q$ is Scott convexity-preserving if and only if $f:(P, \sigma^{\ast}(P))\longrightarrow(Q, \sigma^{\ast}(Q))$ is convexity-preserving.
\end{pn}
\noindent {\bf Proof.}
To check the  necessity, let $B\in \sigma^{\ast}(Q)$. Since $f$ is $L$-order-preserving, we have that $f^{\leftarrow}(B)$ is a lower set. For every nonempty finite $L$-subset $F$ of $P$, $f^{\rightarrow}(F)$ is a nonempty finite $L$-subset of $Q$ since $f^{\rightarrow}\dashv f^{\leftarrow}$. If $F$ has a supremum, we have
\begin{align*}
{\rm sub}(F, f^{\leftarrow}(B))&={\rm sub}(f^{\rightarrow}(F), B)\leq B(\sqcup f^{\rightarrow}(F)) \\
&=B(f(\sqcup F))\\
&=f^{\leftarrow}(B)(\sqcup F).
\end{align*}
Thus, $f^{\leftarrow}(B)\in \sigma^{\ast}(P)$. This shows that $f:(P, \sigma^{\ast}(P))\longrightarrow(Q, \sigma^{\ast}(Q))$ is convexity-preserving.

In what follows, we  prove the sufficiency. For all  $a, b\in P$,  since ${\downarrow}f(b)\in \sigma^{\ast}(Q)$, we have $f^{\leftarrow}({\downarrow}f(b))\in \sigma^{\ast}(P)$. Hence, $f^{\leftarrow}({\downarrow}f(b))$ is a lower set. Thus,  we have
$$e(a, b)=e(a, b)\otimes f^{\leftarrow}({\downarrow}f(b))(b)\leq f^{\leftarrow}({\downarrow}f(b))(a)=e(f(a), f(b)).$$
This show that $f$ is $L$-order-preserving. Let $F$ be a nonempty finite $L$-subset with a supremum. We claim that $f(\sqcup F)=\sqcup f^{\rightarrow}(F)$. In fact, for very $y\in Q$,
$${\rm sub}(f^{\rightarrow}(F), {\downarrow} y)={\rm sub}(F, f^{\leftarrow}({\downarrow} y))
=f^{\leftarrow}({\downarrow} y)(\sqcup F)
=e(f(\sqcup F), y).
$$
Thus, $f(\sqcup F)=\sqcup f^{\rightarrow}(F)$. This shows that $f:P\longrightarrow Q$ is Scott convexity-preserving.
\hfill$\Box$\\

{In \cite{Mao-Xu-2005}, Mao and Xu provided a topological representation of directed completions of consistent algebraic L-domains via sobrification. Recently, Zhang et al. \cite{Zhang-Wang} studied the relationship between  D-completions of  $L$-cotopological space and fuzzy directed completions of $L$-ordered sets. In \cite{Zhang-Wang}, fuzzy directed sets are crucial in defining   monotone convergence $L$-cotopological spaces, also called  fuzzy $d$-spaces. While in this paper, the definition of sobriety of $L$-convex spaces
is introduced using fuzzy nonempty finite sets. These motivate us to  investigate  the  relationship between sober  $L$-convex spaces  and a certain type of fuzzy finite complete $L$-ordered sets.}

\begin{dn}  An $L$-ordered set $P$ is called an {\em $L$-join-semilattice} if every nonempty finite $L$-subset $F$  of $P$ has a supremum.

\end{dn}

\begin{pn}\label{pn1-sober-join} Let $(X, \mathcal{C})$ be a sober $L$-convex space. Then $\Omega X$ is an $L$-join-semilattice.
\end{pn}
\noindent {\bf Proof.} Let $F\in L^{X}$ be a nonempty finite $L$-subset.  Since $X$ is sober, there exists a unique $a\in X$ such that $co(F)=co(1_a)$.  We claim that $a$ is the supremum of $F$ in $\Omega X$. In fact, by Proposition \ref{pn-spec-co}
\begin{align*}
e_{\mathcal{C}(X)}(a, x)&= co(1_x)(a)= {\rm sub}(1_a, co(1_x))\\
&= {\rm sub}(co(1_a), co(1_x))\\
&={\rm sub}(co(F), co(1_x))\\
&={\rm sub}(F, co(1_x))\\
&={\rm sub}(F, {\downarrow} x).
\end{align*}
This shows that $a$ is the supremum of $F$. Thus $\Omega X$ is an $L$-join-semilattice.\hfill$\Box$\\


\begin{pn}\label{pn-sober-join} An $S_0$ $L$-convex space  $(X, \mathcal{C})$ is sober iff $\Omega X$ is an $L$-join-semilattice and $\mathcal{C}\subseteq \sigma^{\ast}(\Omega X)$.
\end{pn}
\noindent {\bf Proof.}
{\bf Necessity}. By {Proposition \ref{pn1-sober-join}}, $\Omega X$ is an $L$-join-semilattice. For every $A\in \mathcal{C}$, let $F$ be a nonempty finite $L$-subset of $X$ ,  we have
$${\rm sub}(F, A)={\rm sub}(co(F), A)={\rm sub}(co(1_{\sqcup F}), A)=A( \sqcup F).$$
This shows that $A\in \sigma^{\ast}(\Omega X)$. Thus $\mathcal{C}\subseteq \sigma^{\ast}(\Omega X)$.

{\bf Sufficiency}. Let $F$ be a nonempty finite $L$-subset of $X$.  Since $(X, e_{\mathcal{C}})$ is an  $L$-join-semilattice, the supremum of $F$ exists.  We claim that $co(F)=co(1_{\sqcup F})$. In fact, since $co(F)\in \mathcal{C}\subseteq\sigma^{\ast}(\Omega X)$ and $F\leq co(F)$, we have $co(F)(\sqcup F)=1$. That is to say, $1_{\sqcup F}\leq co(F)$. Thus $co(1_{\sqcup F})\leq co(F)$. On the other hand, since
$$ 1=e_{\mathcal{C}}(\sqcup F, \sqcup F)=\bigwedge_{x\in X}F(x)\longrightarrow {\downarrow} \sqcup F(x)=\bigwedge_{x\in X}F(x)\to  co(1_{\sqcup F})(x),                $$
we have $F\leq co(1_{\sqcup F})$. Thus $co(F)=co(1_{\sqcup F})$,  which shows that $(X, \mathcal{C})$ is a  sober $L$-convex space.\hfill$\Box$\\

{Note the difference between the above theorem and \cite[Theorem 3.13]{Yao-Zhou}. \cite[Theorem 3.13]{Yao-Zhou} shows that for
 every classical sober convex space $(X, \mathcal{C})$,  $\mathcal{C}=\sigma^{\ast}(\Omega X)$. But in the fuzzy setting,  $\sigma^{\ast}(\Omega X)\subseteq\mathcal{C}$ cannot be obtained.
Moreover, by the proof of Proposition \ref{pn-sober-join}, we know that the Scott $L$-convex structure on an $L$-join-semilattice is indeed sober.}

Unless otherwise stated,  whenever an order-theoretic concept is mentioned in the context of an  $S_0$ $L$-convex space $X$, it is to be interpreted with respect to the specialization $L$-order on $X$. For example, we use  order-theoretic concept to restate Definition \ref{dn-F-close} in the case of sober $L$-convex spaces:
Let $X$ be a sober $L$-convex space and $A\subseteq X$. Then $A$ is called an {\em F-closed set}, if for every nonempty finite $L$-subset with $F\leq \chi_A$ implies $\sqcup F\in A$.

Moreover,  we use  order-theoretic concept to recall Proposition \ref{pn-fcon}(1), {we deduce that:}
 Let $X$ be a sober $L$-convex space and $Y$ be an   $S_0$ $L$-convex spaces.  If $f:X\longrightarrow Y$ be a convexity-preserving mapping,  then $f:\Omega X\longrightarrow \Omega Y$ is Scott convexity-preserving.
\vskip 3pt

We use $(P^F, \mathcal{C}(P^F))$ with the mapping $\xi_P$  to denote the sobrification of $(P, \sigma^{\ast}(P))$.
\begin{lm}\label{lm-xi-scot}
The mapping $\xi_{P}:(P, e_P)\longrightarrow (P^F, {\rm sub}_P)$ is Scott convexity-preserving.
\end{lm}
\noindent {\bf Proof.}
If $K$ is a nonempty finite $L$-subset of $P$ and the supremum $\sqcup K$ exists, then $\sqcup\xi_{P}^{\rightarrow}(K)=co(1_{\sqcup K})(=\xi_{P}(\sqcup K))$. In fact, for every $B\in P^F$,
   \begin{align*}
{\rm sub}(\xi_{P}^{\rightarrow}(K), {\rm sub}(-, B))&=\bigwedge_{x\in X} \xi_{P}^{\rightarrow}(K)(co(1_x))\to  {\rm sub}(co(1_x), B)\\
&=\bigwedge_{x\in X} K(x)\to B(x)\\
&={\rm sub}(K, B)\\
&=B(\sqcup K)\\
&= {\rm sub}(co(1_{\sqcup K}), B).
\end{align*}
This shows that $\xi_{P}$ is Scott convexity-preserving.
\hfill$\Box$\\


By proposition \ref{pn-xf-spe}, when see $P^F$ as an $L$-ordered set, the related $L$-order is always assumed to be the specialization $L$-order of $(P^F, \mathcal{C}(P^F))$; that is inclusion $L$-order ${\rm sub}_P$. The following proposition shows that the original $L$-convex structure of  the sobrification $P^F$ coincides with the Scott $L$-convex structure of  $(P^F, {\rm sub}_P)$.

\begin{pn}\label{pn-c-sig}Let $P$ be an $L$-ordered set. Then $\mathcal{C}(P^F)=\sigma^{\ast}(P^F)$.
\end{pn}
\noindent {\bf Proof.}
By Proposition \ref{pn-sober-join},  we only need to    show $\sigma^{\ast}(P^F)\subseteq\mathcal{C}(P^F)$. It follows from \ref{pn-sober-join} that $(P^F, {\rm sub}_P)$ is an $L$-join-semilattice and $(P, \sigma^{\ast}(P^F))$ is a sober $L$-convex space. By Lemma \ref{lm-xi-scot}, $\xi_P:(P, \sigma^{\ast}(P))\longrightarrow(P^F, \sigma^{\ast}(P^F))$ is convexity-preserving.
Since $(P^F, \mathcal{C}(P^F))$ with mapping $\xi_P$ is a sobrification of $(P, \sigma^{\ast}(P))$, there exists a unique convexity-preserving mapping $\overline{\xi_P}:(P^F, \mathcal{C}(P^F))\longrightarrow (P^F, \sigma^{\ast}(P^F))$ such that the following diagram commutes.
\begin{displaymath}
\xymatrix@=8ex{(P, \sigma^{\ast}(P))\ar[r]^{\xi_P}\ar[dr]_{\xi_P}&(P^F, \mathcal{C}(P^F))\ar@{-->}[d]^{\overline{\xi_P}}\\  &(P^F, \sigma^{\ast}(P^F))}
\end{displaymath}
Write $\mathcal{A}=\{A\in P^F\mid \overline{\xi}_P(A)=A\}$. It is clear that $\Theta(P)\subseteq \mathcal{A}$. We next prove that $\mathcal{A}$ is an F-closed set of
$(P^F, \mathcal{C}(P^F))$. {Here, we need to note that} the specialization $L$-orders of $(P^F, \mathcal{C}(P^F))$ and $(P^F, \sigma^{\ast}(P^F))$ are consistent; that is ${\rm sub}_{P}$. Let $\mathcal{K}$ be a nonempty finite $L$-subset of $ P^F$ and $\mathcal{K}\leq \chi_{\mathcal{A}}$.  By Proposition \ref{pn-fcon}(1), we have $\sqcup \overline{\xi}_{P}^{\rightarrow}(\mathcal{K})= \overline{\xi}_{P}(\sqcup\mathcal{K})$.  Since $\mathcal{K}\leq \chi_{A}$, we have $\overline{\xi}_{P}^{\rightarrow}(\mathcal{K})=\mathcal{K}$. Thus, $\sqcup \overline{\xi}_{P}^{\rightarrow}(\mathcal{K})= \sqcup\mathcal{K}$. Hence, $\overline{\xi}_{P}(\sqcup\mathcal{K})=\sqcup\mathcal{K}$.   This shows that $\sqcup \mathcal{K}\in \mathcal{A}$. Thus $\mathcal{A}$ is an F-closed set of
$(P^F, \mathcal{C}(P^F))$. By Lemma \ref{lm-Fclo-two}, $\mathcal{A}$ is also an F-closed set of
${\rm Cp}(\sigma^{\ast}(P))$. Noticing that $\Theta(P)\subseteq\mathcal{A}\subseteq P^F$, we have $\mathcal{A}= P^F$. Thus $\overline{\xi_p}=id_{P^F}$. Since $\overline{\xi_p}$ is convexity-preserving, we have $\sigma^{\ast}(P^F)\subseteq\mathcal{C}(P^F)$, as desired.
\hfill$\Box$\\
\begin{dn}
 Let $P$ and  $Q$ be two $L$-ordered sets. If $Q$ is an $L$-join-semilattice and $j:P\longrightarrow Q$ is Scott convexity-preserving, then   $(Q, j)$, or $Q$, is called an  {\em $L$-join-semilattice completion} of $P$ if for every  $L$-join-semilattice $M$ and Scott convexity-preserving mapping $f: P\longrightarrow M$, there exists a unique Scott convexity-preserving mapping $\overline{f}: Q\to M$ such that $\overline{f}\circ j=f$.
\end{dn}
By the  universal property of $L$-join-semilattice  completions,  the  $L$-join-semilattice completion of an $L$-ordered set  is unique up to $L$-order-isomorphism.

\begin{tm} \label{tm-join-ord} $(P^{F},{\rm sub}_P)$ with $\xi_{P}$ is an  $L$-join-semilattice completion of $(P,e_p)$.
\end{tm}
\noindent {\bf Proof.}
Let $(Q, e_Q)$ be an  $L$-join-semilattice and Let $f:P\longrightarrow Q$ be Scott convexity-preserving. Then by Proposition \ref{pn-sober-join}$, (Q, \sigma^{\ast}(Q))$ is a  sober $L$-convex space and
$f:(P, \sigma^{\ast}(P))\longrightarrow (Q, \sigma^{\ast}(Q))$ is convexity-preserving. By Theorem \ref{tm-sob-ri} and  Proposition \ref{pn-c-sig},  we have  $(P^F, \sigma^{\ast}(P^F))$  {with $\xi_{P}$} is a sobrification of $(P, \sigma^{\ast}(P))$. Thus,  there exists a unique convexity-preserving mapping $\overline{f}:(P^F, \sigma^{\ast}(P^F))\longrightarrow (Q, \sigma^{\ast}(Q))$ such that $\overline{f}\circ\xi_{P}=f$. By Proposition \ref{pn-sco-dir}, {we have   $\xi_{P}:(P,e_p)\longrightarrow (P^F, {\rm sub}_P)$ and $\overline{f}:(P^F, {\rm sub}_P)\longrightarrow (Q, e_Q)$ are Scott convexity-preserving}. Thus $(P^{F}, {\rm sub}_P)$ with $\xi_{P}$ is an $L$-join-semilattice completion of $(P,e_p)$.
\hfill$\Box$\\

 We now provide a characterization for the  $L$-join-semilattice completion of an $L$-ordered set via sobrification. The proof is straightforward and thus omitted.
 {This result  highlights the significance of our notion of sobriety in bridging fuzzy ordered structures and fuzzy convex structures.}
 \begin{tm} \label{cy-join-ord}
An $L$-ordered set $Q$ is an $L$-join-semilattice completion of an $L$-ordered set $P$ if and only if $(Q, \sigma^{\ast}(Q))$ is a sobrification of $(P, \sigma^{\ast}(P)) $.
\end{tm}
We denote the category of $L$-ordered sets with Scott convexity-preserving mappings as morphisms by $L$-{\bf Ord}.
Theorem \ref{tm-join-ord} shows that the full subcategory of  $L$-join-semilattices is reflective in $L$-{\bf Ord}.

\section{Conclusions}

This paper explore  the connection between sober  $L$-convex spaces and $L$-join-semilattices for a complete residuated lattice $L$ as the truth value table. We construct a  sobrification for any $L$-convex space, demonstrating that the category of sober $L$-convex spaces is a reflective subcategory of that of $L$-convex spaces (cf. Theorem \ref{tm-sob-ri}).   We also show  that   $Q$ is an $L$-join-semilattice completion of $P$ if and only if $(Q, \sigma^{\ast}(Q))$ is a sobrification of $(P, \sigma^{\ast}(P)) $ (cf. Theorems \ref{tm-join-ord},  \ref{cy-join-ord}).
 These work   promote a closer connection between $L$-convex structures  and $L$-ordered structures. Consequently, both ordered and categorical approaches can be effectively combined in future studies of fuzzy convex structure theory.
\vskip 2mm
 We offer two future research directions as follows:
\begin{itemize}

\item {In the classical case, the notions of  polytopes and compact convex sets are equivalent. However,  in the fuzzy setting, this equivalence  is not obvious. This is indeed a question to which we possess no answer. Liu and Yue in \cite{Liu-Yue-2024} have defined a type of sobriety for  $L$-convex spaces based on compact convex sets. In this paper, we introduce another type of sobriety based on polytopes. Future research can  focus on the relationship between these  two types of  sobriety.  The first step in bridging this relationship is to either prove the equivalence between compact convex sets and polytopes or provide a counterexample.}


\item  In the theory of topology, studying generalizations of sobriety is  an interesting topic.
Zhao and Fan in \cite{zhao-fan}
 introduced a weaker form  of sobriety, called bounded sobriety.  Zhang and Wang  \cite{Zhang-Wang}
 further extended bounded sobriety   to the framework of $\mathcal{Q}$-cotopological spaces. Following these work, one can introduce the notion of bounded sobriety in the framework of  $L$-convex structure and subsequently investigate the relationship between bounded  sober $L$-convex spaces and conditional  $L$-join-semilattices.
\end{itemize}

\vspace*{.15in}

\noindent{\bf Acknowledgements.}

This paper is supported by the National Natural Science Foundation of China
(12231007, 12371462), Jiangsu Provincial Innovative and Entrepreneurial
Talent Support Plan (JSSCRC2021521).


\vspace*{.15in}



\begin{thebibliography}{99}
	
%

\bibitem{Category} J. Ad\'{a}mek, H. Herrlich, G.E. Strecker, {\it Abstract and Concrete Categories}, Wiley, New York, 1990.
\bibitem{residuated} R. B\v{e}lohl\'{a}vek, {\it Some properties of residuated lattices}, Czech. Math. J. {\bf 53} (2003) 161--171. https://doi.org/10.1023/a:1022935811257

\bibitem{graph}M. Farber, R.E. Jamison, {\it Convexity in graphs and hypergraphs}, SIAM J.  Alg.  Disc. Math. {\bf 7} (1986) 433--444. https://doi.org/10.1137/0607049

\bibitem{graph1}M. Farber, R.E. Jamison, {\it On local convexity in graphs}, Discrete Math. {\bf 66} (1987) 231--247. https://doi.org/10.1016/0012-365x(87)90099-9

\bibitem{Varlet} S.P. Franklin, {\it Some results on order convexity}, Amer. Math. Monthly  {\bf 69} (1962) 357--359. https://doi.org/10.1080/00029890.1962.11989897


\bibitem{Domain} G. Gierz, K.H. Hofmann, K. Keimel, J.D. Lawson, M. Mislove, D.S. Scott, {\it Continuous Lattices and Domains},  Cambridge University Press, New York, 2003.

\bibitem{Hajek} P. H\'{a}jek, {\it Metamathematics of Fuzzy Logic}, Kluwer Academic Publishers, Dordrecht, 1998.
\bibitem{Soltan}  F. Harary, J. Nieminen, {\it  Convexity in graphs}, J.   Dif{}ferential Geom. {\bf 16} (1981) 185--190. https://doi.org/10.4310/jdg/1214436096

\bibitem{MVTop}U. H\"{o}hle, {\it Many-valued Topology and Its Applications}, Kluwer Academic Publishers, New York, 2001.

\bibitem{Stone}P.T. Johnstone, {\it Stone Spaces}, Cambridge University Press, Cambridge, 1982.
\bibitem{Jankowski} A.W. Jankowski, {\it Some modifications of Scott's theorem on injective spaces}, Stud. Log. {\bf 45}  (1986) 155--166. https://doi.org/10.1007/BF00373271
\bibitem{topology2} H. Komiya, {\it Convexity on a topological space}, Fund. Math.  {\bf 111} (1981) 107--113. https://doi.org/10.4064/fm-111-2-107-113
\bibitem{Liu-Meng-Yue} M.Y. Liu, Y.L. Yue, X.W. Wei, {\it Frame-valued Scott open set monad and its algebras}, Fuzzy Sets Syst. {\bf 460} (2023) 52--71.
    
    https://doi.org/10.1016/j.fss.2022.11.002
\bibitem{Liu-Yue-2024} M.Y. Liu, Y.L. Yue, {\it The reflectivity of the category of stratified $L$-algebraic closure spaces}, Iran. J. Fuzzy Syst. {\bf  21} (2024) 117--12.

https://doi.org/10.22111/IJFS.2024.47190.8314


\bibitem{Algebra}E. Marczewski, {\it Independence in abstract algebras results and problems}, Colloq. Math.   {\bf 14} (1966) 169--188. https://doi.org/10.4064/cm-14-1-169-188
\bibitem{Maruyama} Y. Maruyama, {\it Lattice-valued fuzzy convex geometry}, Comput. Geom. Discrete Math. {\bf 164} (2009) 22--37.

\bibitem{Menger}K. Menger, {\it Untersuchungen \"{u}ber allgemeine metrik}, Math. Ann.  {\bf 100} (1928) 75--163. https://doi.org/10.1007/978-3-7091-6110-4\_20
\bibitem{Mao-Xu-2005} X.X. Mao, L.S. Xu, {\it Representation theorems for directed completions of consistent algebraic L-domains}, Algebr. Univ. {\bf 54} (2005)
    
    https://doi.org/10.1007/s00012-005-1953-x

\bibitem{algebra4}J. Nieminen, {\it The ideal structure of simple ternary algebras}, Colloq. Math. {\bf 40} (1978) 23--29. https://doi.org/10.4064/cm-40-1-23-29



\bibitem{Pang-2023}B. Pang, {\it Quantale-valued convex structures as lax algebras}, Fuzzy Sets Syst. {\bf 473} (2023) 108737.  https://doi.org/10.1016/j.fss.2023.108737
\bibitem{PangShi}B. Pang, F.G. Shi, {\it Subcategories of the category of $L$-convex spaces}, Fuzzy Sets Syst. {\bf 313} (2017) 61--74.  https://doi.org/10.1016/j.fss.2016.02.014

\bibitem{Pang2}B. Pang,  F.G. Shi, {\it Fuzzy counterparts of hull operators and interval operators in the framework of L-convex spaces}, Fuzzy Sets Syst. {\bf 369} (2019) 20--39. https://doi.org/10.1016/j.fss.2018.05.012

\bibitem{Pang3}B. Pang, Z.Y. Xiu, {\it An axiomatic approach to bases and subbases in L-convex spaces and their applications}, Fuzzy Sets Syst. {\bf 369} (2019) 40--56. https://doi.org/10.1016/j.fss.2018.08.002
%


\bibitem{RosaDc} M.V. Rosa, {\it A study of fuzzy convexity with special reference to separation properties}, Cochin University of science and Technology, Cochin, India, 1994.


\bibitem{Quantale}K.I. Rosenthal, {\it Quantales and Their Applications},
Longman House, Burnt Mill, Harlow, 1990.


\bibitem{Shen}C. Shen, S.J. Yang, D.S. Zhao, F.G. Shi, {\it Lattice-equivalence of convex spaces}, Algebra Univer. {\bf 80} (2019)  26.  https://doi.org/10.1007/s00012-019-0600-x

\bibitem{ShenShi} C. Shen, F.G. Shi, {\it $L$-convex systems and the categorical isomorphism to Scott-hull operators},  Iran. J.  Fuzzy Syst. {\bf 15}  (2018) 23--40. 

https://doi.org/10.22111/ijfs.2017.3296



\bibitem{Xiu2} F.G. Shi, Z.Y. Xiu, {\it $(L,M)$-fuzzy convex structures}, J. Nonlinear Sci. Appl. {\bf 10} (2017) 3655--3669.  https://doi.org/10.22436/jnsa.010.07.25



\bibitem{Van1984} M. Van De Vel, {\it Binary convexities and distributive lattices}, Proc.  London  Math. Soc. {\bf 48} (1984) 1--33.  https://doi.org/10.1112/plms/s3-48.1.1

\bibitem{topology1} M. Van De Vel, {\it On the rank of a topological convexity}, Fund. Math.  {\bf 119} (1984) 17--48.  https://doi.org/10.4064/fm-119-2-101-132

\bibitem{Vanbook} M. Van De Vel, {\it Theory of convex spaces}, North-Holland, Amsterdam, 1993.

\bibitem{WangKai}K. Wang, F.G. Shi, {\it Many-valued convex structures induced by
fuzzy inclusion orders}, J. Intell. Fuzzy Syst. {\bf 36} (2019) 2705--2713.  https://doi.org/10.3233/jifs-181103
\bibitem{XU-TopAPpp} L.S. Xu, {\it Continuity of posets via Scott topology and sobrifiction}, Topology Appl. {\bf 153} (2006) 1886--1894. https://doi.org/10.1016/j.topol.2004.02.024
\bibitem{Xia-FFS} C.C. Xia, {\it A categorical isomorphism between injective balanced $L$-$S_0$-convex spaces and fuzzy frames}, Fuzzy Sets Syst. {\bf 437} (2022) 114--126. https://doi.org/10.1016/j.fss.2021.09.018
\bibitem{Xia-further}C.C. Xia, {\it Some further results on pointfree convex geometry}, Algebra Univers. {\bf 85} (2024)
 20. https://doi.org/10.1007/s00012-024-00847-7

\bibitem{YaoTFS}W. Yao, {\it A categorical isomorphism between injective stratified fuzzy $T_0$-spaces and fuzzy continuous lattices}, IEEE Trans. Fuzzy Syst. {\bf 24} (2016) 131--139. https://doi.org/10.1109/tfuzz.2015.2428720
\bibitem{PartI}W. Yao, {\it Quantitative domains via fuzzy sets: Part I: Continuity
of fuzzy directed-complete poset}, Fuzzy Sets Syst. {\bf 161} (2010), 983--987.  https://doi.org/10.1016/j.fss.2009.06.018
\bibitem{PartII}W. Yao, F.G. Shi, {\it Quantitative domains via fuzzy sets:
Part II: Fuzzy Scott topology on fuzzy directed complete posets},
 Fuzzy Sets Syst. {\bf 173} (2011) 60--80. https://doi.org/10.1016/j.fss.2011.02.003

\bibitem{Yao-Yue-FFS}W. Yao, Y.L. Yue,  {\it Algebraic representation of frame-valued continuous lattices via open filter monad}, Fuzzy Sets Syst. {\bf 420} (2021) 143--156.  https://doi.org/10.1016/j.fss.2021.02.004
\bibitem{Yue-Yao-Ho}Y.L. Yue, W. Yao, W.K. Ho, {\it Applications of Scott-closed sets in convex structures}, Topol.  Appl. {\bf 314} (2022) 108093.

    https://doi.org/10.1016/j.topol.2022.108093
\bibitem{Yao-Zhou}W.Yao, C.J. Zhou, {\it Representation of sober convex spaces by join-semilattices}, J. Nonlinear Convex Anal. {\bf 21} (2020) 2715-2724.
\bibitem{Yao-Zhou-2021}W.Yao, C.J. Zhou, {\it A lattice-type duality of lattice-valued fuzzy convex spaces}, J. Nonlinear Convex Anal.
 {\bf 21} (2021) 2843--2853
\bibitem{DX-Zhang-2018} D.X. Zhang, {\it Sobriety of quantale-valued cotopological spaces}, Fuzzy Sets Syst. {\bf 350} (2018) 1--19. https://doi.org/10.1016/j.fss.2017.09.005
\bibitem{FanZhang}Q.Y. Zhang, L. Fan, {\it Continuity in quantitative domains},  Fuzzy Sets Syst. {\bf 154} (2005) 118--131. https://doi.org/10.1016/j.fss.2005.01.007
\bibitem{zhao-fan}D.S. Zhao, T.H. Fan, {\it Dcpo-completion of posets}, Theor. Comput. Sci. {\bf 411} (2010)  2167--2173. https://doi.org/10.1016/j.tcs.2010.02.020
\bibitem{Zhang-Shi-Li} Z.X. Zhang, F.G. Shi, Q.G. Li, K. Wang, {\it On fuzzy monotone convergence $\mathcal{Q}$-cotopological spaces}, Fuzzy Sets Syst. {\bf 425} (2021) 18--33.
    
     https://doi.org/10.1016/j.fss.2020.11.021
\bibitem{Zhang-Wang} Y. Zhang, K. Wang, {\it Bounded sobriety and k-bounded sobriety of $\mathcal{Q}$-cotopological spaces}, Filomat  {\bf 33}  (2019) 2095--2106. https://doi.org/10.2298/fil1907095z


\end{thebibliography}
\end{document}